\documentclass[11pt]{article}
\usepackage{bez123,calc,curves,ebezier,epic,eepic,graphicx,multiply,rotating}
\everymath{\displaystyle}
\usepackage{graphicx}
\usepackage{pst-all}
\usepackage{color}
\usepackage{amssymb}
\usepackage{amsmath}
\newtheorem{prethm}{{\bf Theorem}}

\newenvironment{thm}{\begin{prethm}{\hspace{-0.5
               em}{\bf .}}}{\end{prethm}}
\newtheorem{prelemma}{{\bf Lemma}}

\newtheorem{preex}{{\bf Example}}

\newtheorem{preprop}{{\bf Proposition}}

\newenvironment{prop}{\begin{preprop}{\hspace{-0.5em}{\bf .}}}{\end{preprop}}
\newtheorem{precor}{{\bf Corollary}}

\newtheorem{preremark}{{\bf Remark}}

\newtheorem{preprob}{{\bf Problem}}

\newtheorem{predefin}{{\bf Definition}}

\newtheorem{preconj}{{\bf Conjecture}}

\newtheorem{preprobb}{{\bf Problem}}

\newtheorem{prelem}{{\bf Theorem}}

\newenvironment{proof}{{\bf Proof.}\rm }{\hfill{$\Box$}}

\newtheorem{presolution}{{\bf Solution.}}

\def\newpic#1{}
\def\qed{\ifhmode\unskip\nobreak\fi\quad\ifmmode\Box\else$\Box$\fi}

\title{\vspace{-0.1cm}\Large\bf Generalized degeneracy, dynamic monopolies and maximum degenerate subgraphs}
\author{\large\bf Manouchehr Zaker\footnote{E-mail: mzaker@iasbs.ac.ir}
\vspace{5mm}\\
    Department of Mathematics,\\
     Institute for Advanced Studies in Basic Sciences,\\
     Zanjan 45137-66731, Iran}

    \date{}

\begin{document}
\maketitle\
\begin{abstract}
\noindent A graph $G$ is said to be a $k$-degenerate graph if any subgraph of $G$ contains a vertex of degree at most $k$. Let $\kappa$ be any non-negative function on the vertex set of $G$. We first define a $\kappa$-degenerate graph. Next we give an efficient algorithm to determine whether a graph is $\kappa$-degenerate. We revisit the concept of dynamic monopolies in graphs. The latter notion is used in formulation and analysis of spread of influence such as disease or opinion in social networks. We consider dynamic monopolies with (not necessarily positive) but integral threshold assignments. We obtain a sufficient and necessary relationship between dynamic monopolies and generalized degeneracy. As applications of the previous results we consider the problem of determining the maximum size of $\kappa$-degenerate (or $k$-degenerate) induced subgraphs in any graph. We obtain some upper and lower bounds for the maximum size of any $\kappa$-degenerate induced subgraph in general and regular graphs. All of our bounds are constructive.
\end{abstract}


\noindent {\bf Keywords:} degeneracy, dynamic monopolies.


\section{Introduction}

\noindent All graphs in this paper are finite, undirected and simple. For standard graph theoretical notions and notations we refer the reader to \cite{BM}. Let $k$ be any non-negative integer. A graph $G$ is said to be a $k$-degenerate graph if any subgraph of $G$ contains a vertex of degree at most $k$. It is a well-known fact that $G$ is $k$-degenerate if and only if the vertices of $G$ can be ordered as $v_1, v_2, \ldots, v_n$ such that the degree of $v_i$ in the subgraph of $G$ induced by $\{v_1, \ldots, v_{i}\}$ is at most $k$ for any $i\in \{1, \ldots, n\}$. The concept of degeneracy has many applications in graph theory for example in extremal and chromatic graph theory and has caused many interesting results and problems in these areas (see e.g. \cite{JT}).

\noindent We generalize the concept of $k$-degeneracy as follows. Let $\kappa$ be any assignment of non-negative integers to the vertices of $G$. We say a graph $G$ is $\kappa$-degenerate if the vertices of $G$ can be ordered as $v_1, v_2, \ldots, v_n$ such that the degree of $v_i$ in the subgraph of $G$ induced by $\{v_1, \ldots, v_{i}\}$ is at most $\kappa(v_i)$ for any $i\in \{1, \ldots, n\}$. Note that when $\kappa$ is a constant function then $\kappa$-degeneracy is equivalent to $k$-degeneracy.

\noindent The other concept to be discussed in this paper is dynamic monopoly. In recent years, great attentions have been paid to the modeling and analysis of the spread of belief or influence in social networks. The concept of dynamic monopolies were introduced in order to formulate these problems. By a threshold assignment for a graph $G$ we mean any function
$\tau: V(G)\rightarrow \Bbb{Z}$ such that $\tau(v)\leq deg(v)$ for any vertex $v$, where $V(G)$ is the vertex set of $G$ and $deg(v)$ is the degree of $v$ in $G$. A subset $D\subseteq V(G)$ is called a dynamic monopoly (or $\tau$-dynamic monopoly) if there exists a partition of $V(G)$ into subsets $D_0, D_1, \ldots, D_k$ such that $D_0=D$ and for any $i=1, \ldots, k-1$ each vertex $v$ in $D_{i+1}$ has at least $\tau(v)$ neighbors in $D_0\cup \ldots \cup D_i$. When a vertex $v$ belongs to $D_i$ for some $i$, we say $v$ is activated at time step $i$. Note that we have extended the notion of dynamic monopolies by allowing the vertices of the graph to have non-positive thresholds. In the standard definition of dynamic monopolies all thresholds are non-negative integers. When a vertex $v$ has threshold $\tau(v) \leq 0$ then $v$ is automatically an active vertex and it belongs to the monopoly; but in the definition of the size of a dynamic monopoly such vertices are not counted. In precise words, by the size of a dynamic monopoly $D=D_0, D_1, \ldots, D_k$ we mean $|D\setminus \{v: \tau(v)\leq 0\}|$. Dynamic monopolies have been widely studied by various authors with various types of threshold assignments \cite{ABW, CL, C, DR, FKRRS, KSZ, Z}. For more related works on dynamic monopolies we refer the reader to \cite{Z} where dynamic monopolies with general threshold assignments were introduced. As we mentioned before, we have extended the notion of dynamic monopolies by allowing the vertices of the graph to have non-positive thresholds in order to obtain more applications of dynamic monopolies. One of these applications is the relationship of dynamic monopolies with $\kappa$-degeneracy to be explored in this paper and in finding $k$-degenerate induced subgraphs with maximum cardinality.

\noindent In the following by $(G,\tau)$ we mean any graph $G$ together with an assignment of thresholds $\tau$ for the vertices of $G$. The concept of resistant subgraphs were introduced in \cite{Z} as follows. A subgraph $K$ of $(G,\tau)$ is resistant subgraph if for any vertex
$v\in K$ one has $deg_K(v)\geq deg_G(v)-\tau(v)+1$, where $deg_K(v)$ is the degree of $v$ in $K$. Note that if a vertex $v$ is such that $\tau(v)\leq 0$ then $v$ automatically belongs to any dynamic monopoly. This means that $v$ does not belong to any resistant subgraph. The following proposition interprets the concept of dynamic monopolies in terms of resistant subgraphs. We have the following result from \cite{Z}.

\begin{thm}
A subset $D$ in $(G,\tau)$ is dynamic monopoly if and only if $G\setminus D$ does not contain any resistant subgraph.\label{resist}
\end{thm}

\section{$\kappa$-degeneracy}

Inspired by Theorem \ref{resist} we obtain the following result.

\begin{thm}
Let $\kappa: V(G)\rightarrow \Bbb{N}\cup \{0\}$ and $\tau: V(G)\rightarrow \Bbb{Z}$ be any two functions such that $\tau(v)+\kappa(v)=deg(v)$ for any vertex $v$ of $G$. Then a subset $M$ is $\tau$-dynamic monopoly if and only if $G\setminus M$ is $\kappa$-degenerate.\label{kappatau}
\end{thm}

\noindent \begin{proof}
Assume first that $M$ is a $\tau$-dynamic monopoly and set $H=G\setminus M$. Let the vertices of $H$ be activated according to the following order: $v_h, v_{h-1}, \ldots, v_2, v_1$, where $h=|H|$. We claim that the opposite order of vertices in $H$, i.e. $v_1, v_2, \ldots, v_h$ has the property
$$deg_{H[v_1, \ldots, v_{i-1}]}(v_i)\leq \kappa(v_i).$$
\noindent Since $M$ is dynamic monopoly and $v_i$ is activated after $v_h, \ldots, v_{i+1}$ then
$$deg_{M\cup H[v_h, \ldots, v_{i+1}]}(v_i) \geq \tau(v_i)= deg_G(v_i)-\kappa(v_i).~~~{\bf (1)}$$
\noindent Therefore
$$deg_M(v_i) + deg_{H[v_h, \ldots, v_{i+1}]}(v_i) \geq deg_G(v_i)-\kappa(v_i)$$
\noindent and
$$deg_M(v_i) + deg_H(v_i) - deg_{H[v_1, \ldots, v_{i}]}(v_i) \geq deg_G(v_i)-\kappa(v_i)$$
\noindent so
$$deg_G(v_i) - deg_{H[v_1, \ldots, v_{i}]}(v_i) \geq deg_G(v_i)-\kappa(v_i).$$
\noindent It implies the desired inequality
$$deg_{H[v_1, \ldots, v_{i}]}(v_i)\leq \kappa(v_i).~~~{\bf (2)}$$
\noindent In fact we showed that {\bf (1)} holds if and only if {\bf (2)} holds. Hence the $\kappa$-degeneracy of $H$ implies
that $M$ is $\tau$-dynamic monopoly.
\end{proof}

\noindent It is clear that if $G$ is a $\kappa$-degenerate graph then $|E(G)|\leq \sum_{v\in V(G)} \kappa(v)$. We have the following proposition. In the following we sometimes take the function $\kappa$ equivalent to its range i.e. $\kappa(V(G))$.

\begin{prop}

\noindent (i) Assume that $G$ is $\kappa$-degenerate and $|E(G)| = \sum_{v\in V(G)} \kappa(v)$. Then there exists a vertex $v$ in $G$ such that $deg(v)=\kappa(v)$.

\noindent (ii) Assume that $\kappa$ is given and the vertex $v$ of $G$ is such that $deg(v)=\kappa(v)$. Then $G$ is $\kappa$-degenerate if and only if $G\setminus \{v\}$ is $\kappa'$-degenerate, where $\kappa'$ is obtained by removing $\kappa(v)$ from $\kappa$.\label{prop1}
\end{prop}

\noindent \begin{proof}
To prove (i), since $G$ is $\kappa$-degenerate then there exists $v_1, \ldots, v_n$ such that $deg_{v_1, \ldots, v_i}(v_i)\leq \kappa(v_i)$.
We have $|E(G)|=\sum_i deg_{v_1, \ldots, v_i}(v_i) \leq \sum_{v\in V(G)} \kappa(v) =|E(G)|$. Hence $deg_{v_1, \ldots, v_i}(v_i) = \kappa(v_i)$ for any $i$. In particular $deg_G(v_n)=deg_{v_1, \ldots, v_n}(v_n) = \kappa(v_n)$.

\noindent To prove (ii) we note that if $G\setminus \{v\}$ is $\kappa$-degenerate by the ordering $v_1, v_2, \ldots, v_{n-1}$ then $G$ is $\kappa$-degenerate by the ordering $v_1, \ldots, v_{n-1}, v$. The converse is trivial since if a graph is $\kappa$-degenerate then any subgraph of it is also $\kappa$-degenerate.
\end{proof}

\noindent The following proposition is for the general case $|E(G)| \leq \sum_{v\in V(G)} \kappa(v)$.

\begin{prop}

\noindent (i) Assume that $G$ is $\kappa$-degenerate. Then there exists a vertex $v$ in $G$ such that $0\leq \kappa(v)-deg(v) \leq \sum_{v\in V(G)} \kappa(v) - |E(G)|$.

\noindent (ii) Assume that $\kappa$ is given and the vertex $v$ of $G$ is such that $0\leq \kappa(v)-deg(v) \leq \sum_{v\in V(G)} \kappa(v) - |E(G)|$. Then $G$ is $\kappa$-degenerate if and only if
$G\setminus \{v\}$ is $\kappa'$-degenerate, where $\kappa'$ is obtained by removing $\kappa(v)$ from $\kappa$.\label{prop2}
\end{prop}

\noindent \begin{proof}
To prove (i), since $G$ is $\kappa$-degenerate then there exists $v_1, \ldots, v_n$ such that $deg_{v_1, \ldots, v_i}(v_i)\leq \kappa(v_i)$. We have $|E(G)|=\sum_i deg_{v_1, \ldots, v_i}(v_i) \leq \sum_{v\in V(G)} \kappa(v)$. Obviously $G'=G\setminus \{v_n\}$ is $\kappa'$-degenerate where $\kappa'=\kappa - \kappa(v_n)$. So we write the same inequality for the edges of $G'$. We obtain $|E(G)|-deg_G(v_n)\leq \sum_{i\in \{1, \ldots, n-1\}} \kappa(v_i)$. Hence $v_n$ satisfies the condition of part (i).

\noindent The proof of part (ii) is similar to the proof of part (ii) of Proposition \ref{prop1} and using the argument of part (i) above.
\end{proof}

\noindent Based on Propositions \ref{prop1} and \ref{prop2} we obtain an algorithm which decides if a graph $G$ is $\kappa$-degenerate.

\begin{thm}
There exists an efficient algorithm such that given a graph $G$ and a function $\kappa$, determines whether $G$ is $\kappa$-degenerate.
\end{thm}

\noindent \begin{proof}
While $|E(G)| < \sum_{v\in V(G)} \kappa(v)$ we seek for a vertex $v$ such that $0\leq \kappa(v)-deg(v) \leq \sum_{v\in V(G)} \kappa(v) - |E(G)|$.
If there exists no such vertex then $G$ is not $\kappa$-degenerate by Proposition \ref{prop2}. If there exists such a vertex $v_1$ we replace $G$ by $G_1= G\setminus \{v_1\}$ and replace $\kappa$ by $\kappa_1 = \kappa \setminus \kappa(v_1)$. While $|E(G_i)| < \sum_{v\in V(G_i)} \kappa_i(v)$ we repeat the same procedure and seek for a vertex $v$ with $0\leq \kappa_i(v)-deg_{G_{i}}(v) \leq \sum_{v\in V(G_i)} \kappa_i(v) - |E(G_i)|$. If at some step there exists no such vertex then the algorithm answers ``NO". If there exists such a vertex then we do the same procedure as before.

\noindent Now assume that at some step $i$ of the algorithm we have $|E(G)| = \sum_{v\in V(G)} \kappa(v)$. At this step we seek for a vertex $v$ such that $deg_G(v)=\kappa(v)$. If no such vertex exists then by Proposition \ref{prop1} the algorithm outputs answer ``NO". If there exists such a vertex $v_1$
then we replace $G$ by $G_1= G\setminus \{v\}$ and replace $\kappa$ by $\kappa_1 = \kappa \setminus \kappa(v)$ and repeat the same procedure for $G_1$. By continuing this method if at some step no vertex exists with property $deg_{G_i}(v)=\kappa(v)$ then the algorithm outputs answer ``NO". We repeat this procedure until all vertices of $G$ are checked.
\end{proof}

\section{Maximum $\kappa$-degenerate induced subgraphs}

\noindent The rest of the paper devotes to the applications of Theorem \ref{kappatau} and some results from \cite{Z} and \cite{KSZ} in estimating the size of maximum $\kappa$-degenerate induced subgraphs in general and regular graphs. Finding the maximum induced subgraphs with a given constant degeneracy is an interesting topic in graph theory. The special case of this topic is maximum induced forest or decycling number of graphs. To the best of the author's knowledge, the best known result in this regard is the following theorem of Alon et al. \cite{AKS}. Note that in \cite{AKS}, a graph $G$ is said to be a $k$-degenerate graph if every subgraph of $G$ has a vertex of degree smaller than $k$; but the standard definition of $k$-degenerate graphs is the one we brought in this paper, i.e. we should replace ``smaller than $k$" by ``at most $k$".

\begin{thm} (Alon et al. \cite{AKS})
For every graph $G$, and every integer $k$, there exists a $k$-degenerate induced subgraph $H$ of $G$ such that
$|H|\geq \sum_{v\in V(G)} \min\{1, \frac{k+1}{deg(v)+1}\}$.\label{alon}
\end{thm}


\noindent In order to present our result on maximum $\kappa$-degenerate induced subgraphs, we need the following result from \cite{KSZ} concerning an upper bound for the size of dynamic monopoly with given average threshold. Fortunately the following result still holds even if some threshold values are non-positive.

\begin{thm}
Given any graph $(G,\tau)$ on $n$ vertices and increasing degree sequence $d_1, \ldots, d_n$. There exists an ${\mathcal{O}}(n^3)$ algorithm which outputs a dynamic monopoly $M$ such that $|M|\leq \max \{j: \sum_{i=1}^j (d_i+1) \leq \sum_{v\in V(G)} \
\tau(v)\}$.\label{ksz}
\end{thm}

\noindent Using Theorem \ref{kappatau} and Theorem \ref{ksz} we have the following result.

\begin{thm}
Let $\kappa$ be any positive valued function on the vertex set of $G$. Assume that $G$ has $n$ vertices with increasing degree sequence $d_1, \ldots, d_n$. Then there exists a $\kappa$-degenerate induced subgraph $H$ of $G$ such that
$$|H|\geq n- \max \{j: j+ \sum_{v\in V(G)} \kappa(v) \leq \sum_{i=j+1}^n d_i\}.$$
\noindent Furthermore, such a $\kappa$-degenerate induced subgraph can be obtained by an efficient algorithm.\label{bound}
\end{thm}




\noindent In order to obtain an upper bound for the maximum size of $\kappa$-degenerative induced subgraphs, we use the following result from \cite{Z}.

\begin{thm}
Let $D$ be any $\tau$-dynamic monopoly for $(G,\tau)$. Set $H=G\setminus D$ and let $\tau_{max}$ be the maximum threshold among the vertices of $H$. Then

$\sum_{v\in H}\tau(v) \leq |E(G)|-|E(G[D])|-\delta(G)+t_{max}$.\label{zz}
\end{thm}

\noindent Using Theorem \ref{zz} and Theorem \ref{kappatau} we immediately obtain the following lower bound.

\begin{thm}
Let $\kappa$ be any positive valued function on the vertex set of $G$ with $m$ edges. Set $\alpha= \max\{deg(v)-\kappa(v): v\in V(G)\}$ and $\beta = \min\{deg(v)-\kappa(v): v\in V(G)\}$. Assume that $\beta>0$. Then there exists a $\kappa$-degenerate induced subgraph $H$ such that
 $$|H| \leq \frac{m-\delta(G)+\alpha}{\beta}.$$\label{upperbound}
\end{thm}

\noindent Maximum $k$-degenerate induced subgraphs were studied widely in the literature for various types of graphs such as regular graphs (see e.g. \cite{HW}). In the following we present the results obtained by Theorem \ref{bound} and  Theorem \ref{upperbound} for the size of maximum $k$-degenerate induced subgraph in regular graphs. 

\begin{thm}
Let $G$ be an $r$-regular graph. Denote the size of any $k$-degenerate induced subgraph in $G$ with the maximum cardinality by $\alpha_k(G)$. Then $$\frac{(k+1)n}{r+1} \leq \alpha_k(G)\leq \frac{rn-2k}{2r-2k}.$$\label{regular}
\end{thm}

\noindent The lower bound in Theorem \ref{regular} is the same bound as obtained by applying the bound of Alon et al. for regular graphs. But the advantage of our bounds is that they are all constructive bounds. Theorem \ref{regular} covers a result of Jaeger \cite{J} for cubic graphs asserting that the size of maximum induced forest in any cubic graph is bounded by $(3n-2)/4$.


\end{document}